\documentclass[a4paper,12pt]{article}
\usepackage[top=2.5cm,bottom=2.5cm,left=2.5cm,right=2.5cm]{geometry}
\usepackage{epsfig}
\usepackage{amsmath,amsfonts,amsthm,amssymb}
\usepackage{ifthen}
\usepackage{bbm}  

\newtheorem{theorem}{Theorem}[section]

\newtheorem{corollary}[theorem]{Corollary}

\newcommand{\DEF}[1]{{\em #1\/}}

\newcommand\RR{\ensuremath{\mathbb{R}}}

\newcommand{\R}{\ensuremath{\mathcal{R}}}
\newcommand{\T}{\ensuremath{\mathcal{T}}}

\begin{document}

\title{The Genus Distribution of Doubly Hexagonal Chains}

\author{Bojan Mohar\thanks{Supported in part by an NSERC Discovery Grant, by the
    Canada Research Chair Program, and by the ARRS, Research
    Program P1-0297.}\\[2mm]
    {Department of Mathematics}\\
    {Simon Fraser University}\\
    {Burnaby, B.C., Canada}\\and\\IMFM, Ljubljana, Slovenia\\
    \texttt{mohar@sfu.ca}.}
\date{}
\maketitle

\begin{abstract}
By using the Transfer Matrix Method, explicit formulas for the embedding distribution of doubly hexagonal chain graphs are computed.
\end{abstract}

\begin{center}
Dedicated to Ante Graovac
\end{center}

\baselineskip=0.30in

\section{Fasciagraphs and rotagraphs}

An interesting topic that I have learned from Ante Graovac in the 1980's \cite{BGMP86,BGMP85} is that about fasciagraphs and rotagraphs. Such graphs are frequently studied in crystallography and in mathematical chemistry. They can be described by a small structure that is repeated as long chain that may be either open (fasciagraphs) or closed (rotagraphs). An important tool for dealing with these graphs is a rather general technique based on the \emph{Transfer Matrix Method\/}. In addition to the afore-mentioned articles \cite{BGMP86,BGMP85}, we refer to \cite{Pe69,TeFi88} for some early applications of this method and to \cite{Mu10} for a more recent treatment in theoretical physics.

\begin{figure}[t]
\begin{center}
\includegraphics[width=90mm]{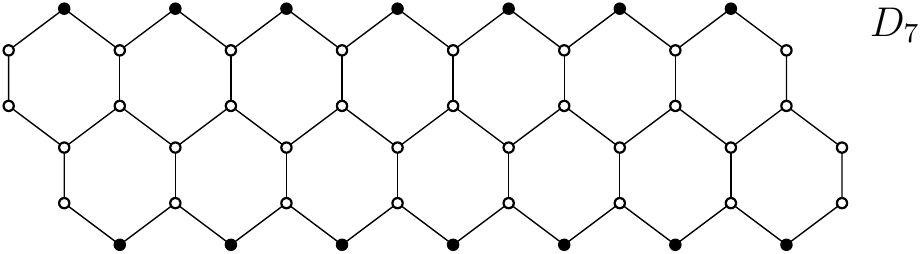}
\end{center}
\caption{The doubly hexagonal chain $D_7$ of length seven}
\label{fig:1}
\end{figure}

Quite recently, the Transfer Matrix Method has been shown to be applicable to the problem of computing genus distributions of fasciagraphs and rotagraphs \cite{Mo13}. In this short note we will illustrate the method on a specific example of fasciagraphs that are isomorphic to the doubly hexagonal chain $D_n$ of length $n$, which is shown in Figure \ref{fig:1} for $n=7$.

We shall restrict ourselves to the following model of fasciagraphs and rotagraphs.
Let $G_0$ and $H$ be fixed graphs, let $S_0\subseteq V(G_0)$ and $S\subseteq V(H)$ be subsets of their vertices, where $|S_0|=|S|$, and let $r:S_0\to S$ be a bijection. Suppose that $B$ is a bipartite graph with bipartition $S_0\cup V(H)$.

We now form a sequence of graphs $G_0,G_1,G_2,\dots$, each with a distinguished vertex set, $S_n\subseteq V(G_n)$ ($n=0,1,2,\dots$), and $S_n$ is in a bijective correspondence with $S_0$. Suppose that we have constructed $G_{n-1}$ and have defined its vertex set $S_{n-1}$, where $S_{n-1}$ is equipped with a bijection to $S_0$. Then we define $G_n$ by first taking the disjoint union of $G_{n-1}$ and a copy of the graph $H$, and then we add the edges of the graph $B$ joining $S_{n-1}$ with $V(H)$, where we identify $S_{n-1}$ with $S_0\subset V(B)$ using the bijection $S_{n-1}\to S_0$. Finally, the vertex-set $S$ in the added copy of the graph $H$ is taken as $S_n$ and $r$ defines its bijection with $S_0$. We say that the sequence $(G_n,S_n)_{n\ge0}$ is a \DEF{linear family} of \DEF{fasciagraphs} with \DEF{constituents} $(G_0,S_0)$, $(H,S)$ and $B$.

\begin{figure}[htb]
\begin{center}
\includegraphics[width=128mm]{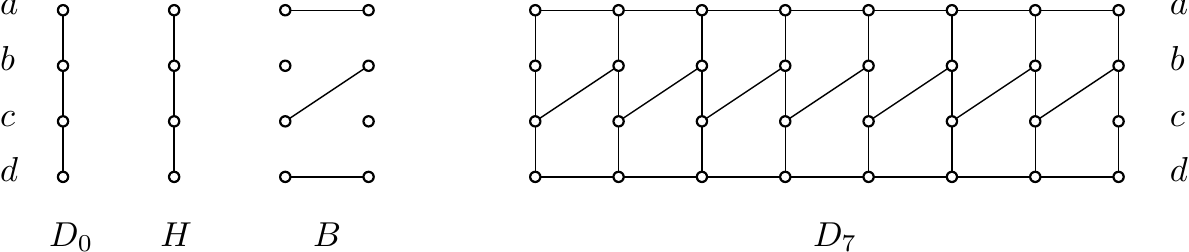}
\end{center}
\caption{Doubly hexagonal chain and its constituents}
\label{fig:2}
\end{figure}

Figure \ref{fig:2} shows the constituents of the doubly hexagonal chain and another drawing of $D_7$.
In this case, we have $D_0=H=P_4$ (the path of length three with vertices $a,b,c,d$), and we have $S_0=S=\{a,b,c,d\}$. For convenience we shall denote the vertices in $S_n$ as $a,b,c,d$ as well.

\begin{figure}[htb]
\begin{center}
\includegraphics[width=52mm]{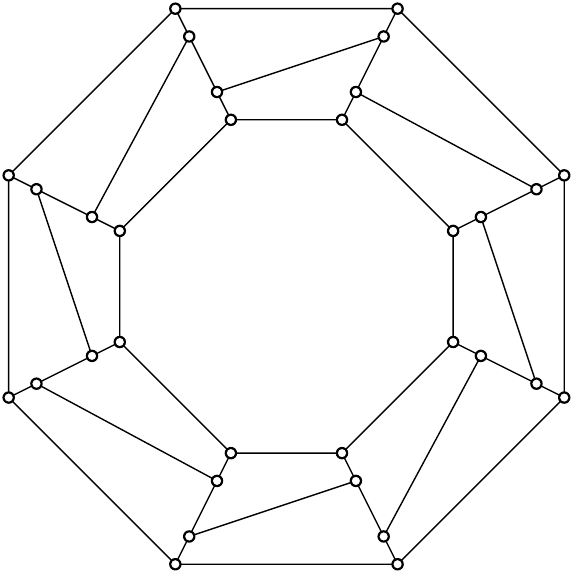}
\end{center}
\caption{Doubly hexagonal rotagraph $R_8$}
\label{fig:3}
\end{figure}

We can identify the last copy of the constituent path $H=P_4$ used to produce $G_n$ with the first copy $D_0$ within $G_n$. In this way we close up the chain and obtain the \emph{rotagraph} $R_n$. Figure \ref{fig:3} shows the graph $R_8$.

\section{Genus distribution}
\label{sect:1}

We speak of a genus distribution of a graph when we consider the collection of all combinatorial types of 2-cell embeddings of the graph in closed orientable surfaces. Due to a variety of applications within mathematics, theoretical physics and theoretical computer science (see, e.g., \cite{LaZv04}), there is abundant literature on the genus distribution of graphs.

Calculations of genus distributions is intimately related to the enumeration of rooted maps on surfaces \cite{MaLiuWe06}. We refer to \cite{LaZv04} for extensive treatment of map enumeration and its applications. Even the very special case of one-vertex graphs (bouquet of loops) bears interesting applications. The one-vertex maps are isomorphic, via surface duality, to unicellular maps that have received certain attention, see, e.g. \cite{BeCh11,Cha10,Cha11}. In summary, embeddings of graphs into surfaces take part in many applications (see, e.g., \cite{LaZv04,MoTh01}); thus embedding distributions are of continuing interest.

Embedding distributions for specific families have been the subject of extensive research. These include different kinds of fasciagraphs (cobblestone path \cite{FuGrSt89,Sta91}, closed-end ladders \cite{Te97}, Ringel ladders \cite{Te00,ChOuZou11}, circular ladders \cite{Li00}, M\"obius ladders \cite{Li05}, lantern graphs \cite{LiHaZh10}, and others \cite{ChLiuWa06,MaLiuWe06,HaHeLiuWe07,YaLiu07,WaLiu08})
as well as rotagraphs (generalized necklaces and circulant necklaces \cite{ShLiu08}, double pearl-ladder graphs \cite{ZeLiu11}).

The 2-cell embeddings of a graph into orientable surfaces can be described combinatorially by means of rotation systems. We refer to \cite{MoTh01} for a thorough treatment and state only the basic definitions here.

A \DEF{rotation system} for a graph $G$ is a set of cyclic permutations, one for each vertex of the graph, $\Pi = \{\pi_v\mid v\in V(G) \}$, where each $\pi_v$ cyclically permutes the edges incident with the vertex $v$. Each rotation system determines a 2-cell embedding of $G$ into an orientable surface so that the clockwise order of the edges emanating from any vertex $v$ matches $\pi_v$. Given $\Pi$, we first determine the set of \DEF{$\Pi$-facial walks} (corresponding to traversals of the boundaries of the faces of the embedding, see \cite{MoTh01} for details), and then we paste a disk along each facial walk in order to obtain a closed surface. The rotation systems on $G$ describe all 2-cell embeddings of $G$ into orientable surfaces up to orientation-preserving homeomorphisms \cite[Theorem 3.2.4]{MoTh01}. Therefore it is customary to speak about (\emph{combinatorial}) \emph{embeddings} instead of rotation systems.

Clearly, the number of cyclic permutations of $\deg(v)$ edges incident with the vertex $v$ is equal to $(\deg(v)-1)!$, and thus the total number of embeddings described by rotation systems is equal to $\prod_{v\in V(G)} (\deg(v)-1)!$. The fasciagraph $D_n$ has $4n-2$ vertices of degree 3, thus it has totally $2^{4n-2}$ combinatorial embeddings.

Let us consider a graph $G$ and its embeddings. For each integer $g\ge0$, let $\alpha_g=\alpha_g(G)$ be the number of (combinatorial) embeddings of $G$ having genus equal to $g$. Let us recall that the \emph{genus} is defined via Euler's formula,
\begin{equation}\label{eq:genus}
  g = 1 + \tfrac{1}{2}(|E(G)| - |V(G)| - f)
\end{equation}
where $f$ is the number of facial walks of the embedding.
The \DEF{genus polynomial} of $G$ is defined as the generating function for the numbers $\alpha_g$, i.e.,
\begin{equation}\label{eq:genus poly}
   \alpha(G,X) = \sum_{g\ge0} \alpha_g\, X^g.
\end{equation}

\section{Embedding types for fasciagraphs $(D_n,S_n)$}
\label{sect:3}

Let $S=S_n\subset V(D_n)$ be the set of vertices of $D_n$ in the last added copy of the constituent $H=P_4$. If $\Pi$ is an embedding of $D_n$ and $W$ is a $\Pi$-facial walk containing at least one vertex of $S$, then the \DEF{$S$-type} of $W$ is the cyclic sequence of vertices of $S$ as they appear in $W$. (We assume for all facial walks that they are given in the cyclic order induced by the orientation of the surface.) We also define the \DEF{$S$-type} of the embedding $\Pi$. This is the family of all $S$-types of $\Pi$-facial walks containing vertices in $S$. Each $S$-type $\tau$ has its reverse $S$-type $\tau^{-1}$ that is obtained by replacing each cyclic sequence in $\tau$ by its inverse. We shall consider the type $\tau$ and its inverse to be the same.

In our case of doubly hexagonal chains, the facial walks containing vertices in $S=\{a,b,c,d\}$ are precisely those facial walks that contain the vertex $b\in S$. These may be three distinct faces, two faces or a single face. If there are three faces, their type is $$\tau_0 = (ba)(abcd)(dcb)$$
(which is considered the same as the inverse $(ab)(dcba)(bcd)$). All together there are six types, the remaining five being:
\begin{align*}
  \tau_1 & = (baabcd)(dcb) \\
  \tau_2 & = (badcb)(abcd) \\
  \tau_3 & = (ba)(abcddcb) \\
  \tau_4 & = (baabcddcb) \\
  \tau_5 & = (badcbabcd)
\end{align*}

\begin{figure}
\begin{center}
\includegraphics[width=94mm]{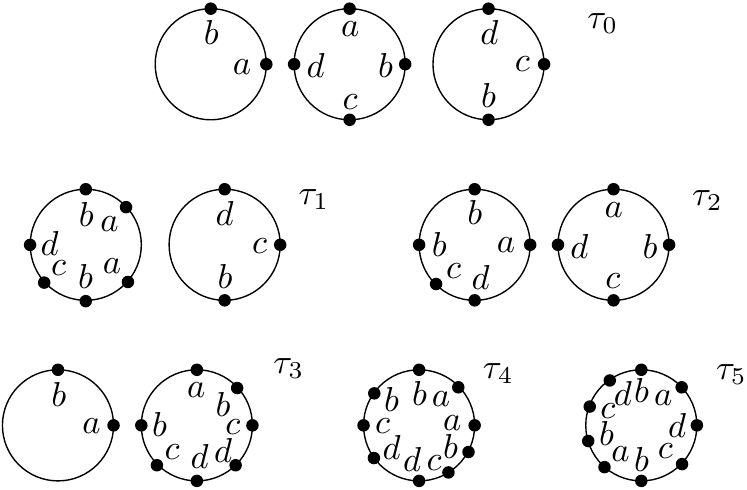}
\end{center}
\caption{The $S$-types for embeddings of $D_n$}
\label{fig:4}
\end{figure}

In Figure \ref{fig:4} we depict a graphical representation of the $S$-types $\tau_i$ ($0\le i\le 5$) by showing the faces containing the vertices in $S$ as disks in the plane.
Observe that each of the vertices $a,c,d\in S$ appears twice (as these vertices have degree 2 in $D_n$) and $b$ appears three times in every $S$-type.

In order to deal with enumeration of embeddings of particular $S$-types, we have to introduce some algebraic notation. We denote by $\T=\{\tau_i\mid 0\le i\le 5\}$ the set of all $S$-types of embeddings of $D_n$ and by $\R$ the set of all functions from $\T$ to $\RR$, $\R = \{\alpha:\T\to\RR\}$.
We consider $\R$ as a 6-dimensional real vector space and also call its elements ``vectors indexed by the $S$-types" and use the vector notation, where the $i$th entry of a vector corresponds to the value $\alpha(\tau_{i-1})$, $1\le i\le 6$.

We can view all embeddings of genus $g$ and consider their types. Then we view $\alpha_g$ as an element of $\R$, by setting $\alpha_g(\tau)$ to be the number of embeddings of genus $g$ whose $S$-type is $\tau\in \T$. To count over different genera we introduce the corresponding formal variable $Y$ and set $$\alpha(D_n,Y,S)=\sum_{g\ge0}\alpha_gY^g.$$
We call $\alpha$ the \DEF{stratified genus distribution}. The difference from (\ref{eq:genus poly}) is that here $\alpha_g$ is no longer a single integer but the stratified genus distribution vector.

\section{Transfer matrices}
\label{sect:join two graphs}

The following theorem has been proved in \cite{Mo13}.

\begin{theorem}
\label{thm:join two graphs}
For every $n\ge2$, every $S$-type $\sigma$ for embeddings of $D_{n-1}$ and for every $S$-type $\tau$ for $D_n$, there is a polynomial $\beta_{\sigma,\tau}(Y)$ which is independent of $n$ such that for every embedding $\Pi$ of $D_{n-1}$ of type $\sigma$ and genus $g$, the number of embeddings of $S_n$ of genus $g+k$, whose induced embedding on $S_{n-1}$ coincides with $\Pi$, is equal to the coefficient of\/ $Y^k$ in $\beta_{\sigma,\tau}(Y)$.
\end{theorem}

Theorem \ref{thm:join two graphs} can be expressed in matrix notation by using the $6\times 6$ matrix $Q$ that is indexed by the set $\T$ of the $S$-types and whose $(\sigma,\tau)$-entry is equal to $\beta_{\sigma,\tau}(Y)$.

\begin{corollary}
\label{cor:path-like}
There is a $6\times 6$ matrix $Q$ such that for every $n\ge2$, we have
$$
   \alpha(D_n,Y,S) = Q^{n-1}\,\alpha(D_1,Y,S).
$$
\end{corollary}

The matrix $Q$ appearing in Corollary \ref{cor:path-like} is called the \DEF{transfer matrix} for the linear family of fasciagraphs $D_n$. If it can be written as, $Q=X^{-1}\Lambda X$, then the form
\begin{equation}\label{eq:diagonalized solution}
    \alpha(D_n,Y,S) = X^{-1}\Lambda^{n-1} X \,\alpha(D_1,Y,S).
\end{equation}
In particular, if $\Lambda$ is a diagonal matrix or is in the Jordan form for $Q$, then (\ref{eq:diagonalized solution}) provides an explicit solution to the recurrence of Corollary \ref{cor:path-like}.

In our example, $D_1$ has two embeddings of genus 0. They are of type $\tau_0$. It also has two embeddings of type $\tau_5$ that are both of genus 1. Thus,
$$
   \alpha(D_1,Y,S) = \left[\begin{array}{c}2\\0\\0\\0\\0\\2Y\end{array}\right]
$$
where the indeterminate $Y=Y^1$ in the last coordinate indicates that the genus is equal to 1.

In order to apply Corollary \ref{cor:path-like}, we need to determine the entries of the transfer matrix $Q$.
For $n\ge2$, going from $D_{n-1}$ to $D_n$ amounts to the same as adding a new vertex $b\in S_n$ and adding three paths from this vertex to the vertices $a,d,c$ in $S_{n-1}$. Thus, for each embedding of $D_{n-1}$ of $S$-type $\tau_i\in\T$, we have 16 possible ways to extend this embedding to an embedding of $D_n$. There are 8 ways, where to put the edges leading from $b$ to the vertices $a,c,d$ (since each of these occurs twice in the facial walks), and then there are two choices for the local rotation around the vertex $b$. If $\tau_i$ has $t\in\{1,2,3\}$ faces, then the resulting embeddings have either $t+2$, $t$ or $t-2$ faces. (The latter possibility occurs only when $t=3$.) Euler's formula (\ref{eq:genus}) shows that if the number of faces increases, then the genus stays the same, if the number of faces stays the same, then the genus increases by 1, and when the number of faces drops from 3 to 1, the genus increases by 2.

From an embedding of type $\tau_4$ or $\tau_5$ we can make three new faces out of one. This occurs in 8 out of 16 possibilities and in these cases the genus stays the same. In the other eight cases, we always obtain an embedding of type $\tau_5$ and the genus increases by one. A similar procedure can be carried over for other types, and this gives us the following transfer matrix
\begin{equation*}
Q =
\left(
  \begin{array}{cccccc}
    1 & 2 & 2 & 4 & 8 & 8  \\
    6Y & 4Y & 4Y & 8Y & 0 & 0  \\
    2Y & 4Y & 4Y & 0 & 0 & 0  \\
    2Y & 4Y & 4Y & 0 & 0 & 0  \\
    2Y^2 & 0 & 0 & 0 & 0 & 0  \\
    2Y^2+Y & 2Y & 2Y & 4Y & 8Y & 8Y  \\
  \end{array}
\right)
\end{equation*}
By Corollary \ref{cor:path-like}, this gives
\begin{equation*}
    \alpha(D_n,Y,S) =
\left(
  \begin{array}{cccccc}
    1 & 2 & 2 & 4 & 8 & 8  \\
    6Y & 4Y & 4Y & 8Y & 0 & 0  \\
    2Y & 4Y & 4Y & 0 & 0 & 0  \\
    2Y & 4Y & 4Y & 0 & 0 & 0  \\
    2Y^2 & 0 & 0 & 0 & 0 & 0  \\
    2Y^2+Y & 2Y & 2Y & 4Y & 8Y & 8Y \\
  \end{array}
\right)^{n-1}
\left[\begin{array}{c}2\\0\\0\\0\\0\\2Y\end{array}\right]
\end{equation*}

From the transfer matrix $Q$ we see that embeddings of types $\tau_1$ and $\tau_2$ behave similarly. Thus we may replace these two $S$-types by a combined type $\tau_{12}$ (where an embedding is said to be of type $\tau_{12}$ if its $S$-type is either $\tau_1$ or $\tau_2$). Similarly, we can replace $\tau_4$ and $\tau_5$ by a combined type $\tau_{45}$. Then the transfer matrix becomes smaller and we obtain the following solution:
\begin{equation}\label{eq:smaller}
    \widehat\alpha(D_n,Y,S) =
\left(
  \begin{array}{cccc}
    1 & 2 & 4 & 8  \\
    8Y & 8Y & 8Y & 0  \\
    2Y & 4Y & 0 & 0  \\
    4Y^2+Y & 2Y & 4Y & 8Y \\
  \end{array}
\right)^{n-1}
\left[\begin{array}{c}2\\0\\0\\2Y\end{array}\right]
\end{equation}

Genus polynomials of $D_n$ for small values of $n$ can be computed directly from (\ref{eq:smaller}) by summing up the four polynomials in the vector $\widehat \alpha(D_n,Y,S)$. They are collected in Table \ref{tbl:1}.

\begin{table}\label{tbl:1}
  \centering
  \begin{tabular}{|c|l|}
    \hline
    $n$ & $\alpha(D_n,x)$ \\
    \hline\hline
    1 & $2+2x$ \\
    2 & $2+38 x+24 x^2$ \\
    3 & $2+102 x+664 x^2+256 x^3$ \\
    4 & $2+166 x+3032 x^2+10368 x^3+2816 x^4$ \\
    5 & $2+230 x+7448 x^2+70912 x^3+152832 x^4+30720 x^5$ \\
    6 & $2+294 x+13912 x^2+244096 x^3+1441536 x^4+2158592 x^5+335872 x^6$ \\[1mm]
    7 & $2+358 x+22424 x^2+595456 x^3+6588672 x^4+26675200 x^5$ \\
    & $~\,+29556736 x^6+3670016 x^7$ \\[1mm]
    8 & $2+422 x+32984 x^2+1190528 x^3+20378368 x^4+155713536 x^5$ \\
    & $~\,+461266944 x^6+395051008 x^7+40108032 x^8$ \\
    \hline
  \end{tabular}
  \caption{Genus polynomials $\alpha(D_n,x)$ for $n\le 8$}\label{tbl:1}
\end{table}

However, the transfer matrix formula also gives way for an explicit solution of the recurrence.
The reduced transfer matrix $\widehat Q$ has the following characteristic polynomial:
$$
  \phi(\widehat Q,\lambda) = \lambda^4-(16Y+1)\lambda^3-16Y\lambda^2+(512Y^3-64Y^2)\lambda+1024 Y^4
$$
The eigenvalues of $\widehat Q$ are:\footnote{The calculation was done with the software package Maple, Version 17.}
\begin{align*}
  \lambda_1 & = 4Y+\frac{1}{4}+\frac{\sqrt{3}}{12}\, w +
  \frac{1}{12} \sqrt{\frac{g+h}{wz}} \\
  \lambda_3 & = 4Y+\frac{1}{4}+\frac{\sqrt{3}}{12}\, w -
  \frac{1}{12} \sqrt{\frac{g+h}{wz}} \\
  \lambda_3 & = 4Y+\frac{1}{4}-\frac{\sqrt{3}}{12}\, w +
  \frac{1}{12} \sqrt{\frac{g-h}{wz}} \\
  \lambda_4 & = 4Y+\frac{1}{4}-\frac{\sqrt{3}}{12}\, w -
  \frac{1}{12} \sqrt{\frac{g-h}{wz}}
\end{align*}
where
\begin{align*}
z &= \Bigl(13824 Y^6 - 864 Y^5 + 207 Y^4 + Y^3 + 9\sqrt{55296 Y^{10} - 2880 Y^9 + 465 Y^8 + 6Y^7}\,\Bigr)^{1/3} \\
w &= \sqrt{(18432 Y^4 + 768 z Y^2 - 768 Y^3 + 32 z^2 + 224 Y z + 32 Y^2 + 3z)/z} \\
g &= -55296 w Y^4 + 4608 w z Y^2 + 2304 w Y^3 - 96 w z^2 + 1344 Y z w - 96 w Y^2 + 18wz \\
h &= 41472 \sqrt{3}\, z Y^2 + 2016 \sqrt{3} \, z Y + 18 \sqrt{3}\, z .
\end{align*}
This implies that
\begin{equation}\label{eq:solution}
    \widehat\alpha(D_n,Y,S) =
N\cdot \left(
  \begin{array}{cccc}
    \lambda_1^{n-1} & 0 & 0 & 0  \\
    0 & \lambda_2^{n-1} & 0 & 0  \\
    0 & 0 & \lambda_3^{n-1} & 0  \\
    0 & 0 & 0 & \lambda_4^{n-1} \\
  \end{array}
\right) \cdot N^{-1}
\left[\begin{array}{c}2\\0\\0\\2Y\end{array}\right]
\end{equation}
where $N$ is the transition matrix for the diagonalization of $\widehat Q$.
This matrix was computed by using software platforms Maple and Mathematica and it turns that
\begin{equation*}
    N =
\left(
  \begin{array}{cccc}
    \varphi_1(\lambda_1) & \varphi_1(\lambda_2) & \varphi_1(\lambda_3) & \varphi_1(\lambda_4) \\[1mm]
    \varphi_2(\lambda_1) & \varphi_2(\lambda_2) & \varphi_2(\lambda_3) & \varphi_2(\lambda_4) \\[1mm]
    \varphi_3(\lambda_1) & \varphi_3(\lambda_2) & \varphi_3(\lambda_3) & \varphi_3(\lambda_4) \\[1mm]
    1 & 1 & 1 & 1 \\
  \end{array}
\right)
\end{equation*}
where $\varphi_1(\lambda) = \frac{\lambda}{Y(4Y+\lambda)}$,
$\varphi_2(\lambda) = \frac{\lambda(\lambda^2 - (8Y+1)\lambda - 32Y^2 -8Y)}{2Y(\lambda^2 + 12Y\lambda + 32Y^2)}$, and
$\varphi_3(\lambda) = \frac{2(\lambda^2-8Y\lambda-32Y^2)}{\lambda^2+12Y\lambda+32Y^2}$.

Note that the $i$th column of $N$ is an eigenvector of $\lambda_i$ ($1\le i\le 4$) that has an alternative representation:
$$
 {\bf u}_i =
 \left[\begin{array}{c}
   -2(\lambda_i^3-24 \lambda_i^2 Y+96 \lambda_i Y^2+512 Y^3-\lambda_i^2-8 \lambda_i Y)\\[1mm]
   \lambda_i^3-24 \lambda_i^2 Y+192 \lambda_i Y^2+512 Y^3-\lambda_i^2-8 \lambda_i Y\\[1mm]
   384Y^3\\[1mm]
   2Y(\lambda_i^3-12 \lambda_i^2 Y+128 Y^3-\lambda_i^2-20 \lambda_i Y-96 Y^2)
 \end{array}\right]
$$
If ${\bf z}_i = {\bf u}_i/\Vert {\bf u}_i\Vert$, then we can write $Q^{n-1}$ in the form
$$
   Q^{n-1} = \sum_{i=1}^4 \lambda_i^{n-1} {\bf z}_i\, {\bf z}_i^T.
$$
This gives rise to an explicit solution for $\widehat\alpha(D_n,Y,S)$ that is equivalent to the one given by (\ref{eq:solution}).

\subsection*{Acknowledgement.}
The author is grateful to Marko Petkov\v sek and to Michael Monagan for their insight and advice on using software packages Mathematica and Maple.

\bibliographystyle{plain}

\bibliography{GenusDistribution_TransferMatrix_Graovac}

\end{document}